\documentclass[11pt,leqno]{article}
\usepackage[doc]{optional} 
\usepackage{xcolor}
\definecolor{labelkey}{rgb}{0,0.08,0.45}
\definecolor{refkey}{rgb}{0,0.6,0.0}
\definecolor{Brown}{rgb}{0.45,0.0,0.05}
\usepackage{mathpazo}
\usepackage{exscale,relsize}
\usepackage{amsmath}
\usepackage{amsfonts}
\usepackage{amssymb}
\usepackage{calc}
\usepackage{theorem}
\usepackage{pifont}      
\usepackage{graphicx}
\DeclareMathOperator{\weakstarly}{\stackrel{\mathrm{w*}}{\rightharpoondown}}

\newcommand{\zeroun}{\ensuremath{\left]0,1\right[}}

\newcommand{\scal}[2]{\langle{{#1},{#2}}\rangle}

\newcommand{\RR}{\ensuremath{\mathbb R}}

\newcommand{\RX}{\ensuremath{\,\left]-\infty,+\infty\right]}}

\newcommand{\NN}{\ensuremath{\mathbb N}}
\newcommand{\nnn}{\ensuremath{{n \in \NN}}}
\newcommand{\thalb}{\ensuremath{\tfrac{1}{2}}}

\newcommand{\menge}[2]{\big\{{#1} \mid {#2}\big\}}

\newcommand{\To}{\ensuremath{\rightrightarrows}}

\newcommand{\spand}{\operatorname{span}}

\newcommand{\dom}{\ensuremath{\operatorname{dom}}}

\newcommand{\gra}{\ensuremath{\operatorname{gra}}}

\newcommand{\intdom}{\ensuremath{\operatorname{int}\operatorname{dom}}\,}
\newcommand{\inte}{\ensuremath{\operatorname{int}}}

\newcommand{\bd}{\ensuremath{\operatorname{bdry}}}

\newcommand{\Id}{\ensuremath{\operatorname{Id}}}

\newcommand{\pinf}{\ensuremath{+\infty}}

\renewcommand{\phi}{\ensuremath{\varphi}}

\newtheorem{theorem}{Theorem}[section]
\newtheorem{lemma}[theorem]{Lemma}
\newtheorem{fact}[theorem]{Fact}
\newtheorem{corollary}[theorem]{Corollary}

\theoremstyle{plain}{\theorembodyfont{\rmfamily}
}
\theoremstyle{plain}{\theorembodyfont{\rmfamily}
}
\theoremstyle{plain}{\theorembodyfont{\rmfamily}
}
\theoremstyle{plain}{\theorembodyfont{\rmfamily}
}
\theoremstyle{plain}{\theorembodyfont{\rmfamily}
\newtheorem{remark}[theorem]{Remark}}
\theoremstyle{plain}{\theorembodyfont{\rmfamily}
}

\begin{document}


\title{\textsc{An Answer to S.\ Simons' Question
on the Maximal Monotonicity of the Sum of a Maximal Monotone
Linear Operator and a Normal Cone Operator }}

\author{
Heinz H.\ Bauschke\thanks{Mathematics, Irving K.\ Barber School,
UBC Okanagan, Kelowna, British Columbia V1V 1V7, Canada. E-mail:
\texttt{heinz.bauschke@ubc.ca}.},\;  Xianfu
Wang\thanks{Mathematics, Irving K.\ Barber School, UBC Okanagan,
Kelowna, British Columbia V1V 1V7, Canada. E-mail:
\texttt{shawn.wang@ubc.ca}.},\; and Liangjin\
Yao\thanks{Mathematics, Irving K.\ Barber School, UBC Okanagan,
Kelowna, British Columbia V1V 1V7, Canada.
E-mail:  \texttt{ljinyao@interchange.ubc.ca}.}}
 \vskip 3mm

\date{February 4, 2009}
\maketitle

\begin{abstract} \noindent
The question whether or not the sum of two maximal monotone
operators is maximal monotone under Rockafellar's constraint
qualification --- that is, whether or not ``the sum theorem'' is true ---
is the most famous open problem in Monotone Operator Theory.
In his 2008 monograph \emph{``From Hahn-Banach to Monotonicity''},
Stephen Simons asked whether or not
the sum theorem holds for the special case of 
a maximal monotone linear operator and a normal cone operator
of a closed convex set provided
that the interior of the set makes a nonempty intersection
with the domain of the linear operator.

In this note, we provide an affirmative answer to Simons' question.
In fact, we show that the sum theorem is true
for a maximal monotone \emph{linear relation} and a normal cone operator.
The proof relies on Rockafellar's formula for the Fenchel conjugate of the
sum as well as some results featuring the Fitzpatrick function. 
\end{abstract}



\section{Introduction}

Throughout this paper, we assume that
$X$ is a Banach space with norm $\|\cdot\|$,
that $X^*$ is its continuous dual space with norm
$\|\cdot\|_*$, and that $\scal{\cdot}{\cdot}$ denotes
the pairing between these spaces. 
Let $A\colon X\To X^*$ 
be a \emph{set-valued operator} (also known as multifunction)
from $X$ to $X^*$, i.e., for every $x\in X$, $Ax\subseteq X^*$,
and let 
$\gra A = \menge{(x,x^*)\in X\times X^*}{x^*\in Ax}$ be
the \emph{graph} of $A$. 
Then $A$ is said to be \emph{monotone} if 
\begin{equation}
\big(\forall (x,x^*)\in \gra A\big)\big(\forall (y,y^*)\in\gra
A\big) \quad \scal{x-y}{x^*-y^*}\geq 0,
\end{equation}
and \emph{maximal monotone} if
no proper enlargement
(in the sense of graph inclusion) of $A$ is monotone.
Monotone operators have proven to be a key class of objects
in modern Optimization and Analysis; see, e.g.,
the books
\cite{BurIus,ph,Si,Si2,RockWets,Zalinescu} and the references therein.
(We also adopt standard notation used in these books:
$\dom A = \menge{x\in X}{Ax\neq\varnothing}$ is the \emph{domain} of $A$.
Given a subset $C$ of $X$, 
$\inte C$ is the \emph{interior},
$\overline{C}$ is the \emph{closure}, 
$\bd{C}$ is \emph{boundary}, and
$\spand C$ is the \emph{span} (the set of all finite linear
combinations) of $C$.
The \emph{indicator function} $\iota_C$ of $C$ takes the value
$0$ on $C$, and $\pinf$ on $X\smallsetminus C$. 
Given $f\colon X\to \RX$, 
$\dom f = f^{-1}(\RR)$ and 
$f^*\colon X^*\to\RX\colon x^*\mapsto\sup_{x\in X}(\scal{x}{x^*}-f(x))$
is the \emph{Fenchel conjugate} of $f$.
Furthermore, $B_X$ is the \emph{closed unit ball}
$\menge{x\in X}{\|x\|\leq 1}$ of $X$, and $\NN=\{0,1,2,3,\ldots\}$.)

Now assume that $A$ is maximal monotone, and 
let $B\colon X\To X^*$ be maximal monotone as well.
While the \emph{sum operator} $A+B\colon X\To X^*\colon x\mapsto
Ax+Bx = \menge{a^*+b^*}{a^*\in Ax\;\text{and}\;b^*\in Bx}$
is clearly monotone, 
it may fail to be maximal monotone. When $X$ is reflexive,
the classical \emph{constraint qualification}
$\dom A \cap\intdom B\neq \varnothing$ guarantees maximal monotonicity
of $A+B$, this is a famous result due to Rockafellar
\cite[Theorem~1]{Rock70}. Various extensions of this \emph{sum theorem} 
have been
found, but the general version in nonreflexive Banach spaces remains
elusive --- this has led to the famous \emph{sum problem}; 
see Simons' recent monograph 
\cite{Si2} for the state-of-the-art. 

The notorious difficulty of the sum problem makes
it tempting to consider various special cases. 
In this paper, we shall focus on the case when
$A$ is a \emph{linear relation} 
and $B$ is the \emph{normal cone operator}
$N_C$ of some nonempty closed convex subset $C$ of $X$.
(Recall that $A$ is a linear relation if $\gra A$ is a linear subspace of
$X\times X^*$, and that for every $x\in X$, the normal cone operator at $x$
is defined by $N_C(x) = \menge{x^*\in
X^*}{\sup\scal{C-x}{x^*}\leq 0}$, if $x\in C$; and $N_C(x)=\varnothing$, 
if $x\notin C$. Consult \cite{Cross} for further information on linear
relations.)
If 
$A\colon X\To X^*$ is \emph{at most single-valued} 
(i.e., for every $x\in X$,
either $Ax=\varnothing$ or $Ax$ is a singleton), then
we follow the common slight abuse of notation to identify $A$ with
a classical operator $\dom A \to X^*$. 
We thus include the classical case when 
$A\colon X\to X^*$ is a continuous linear monotone 
(thus \emph{positive}) operator. 
Continuous and discontinuous linear operators --- and lately even
linear relations --- have received some attention in Monotone
Operator Theory \cite{BB,BBW,BWY2,BWY3,PheSim,Svaiter,VZ} because they provide additional classes of 
examples apart from
the well known and well understood 
\emph{subdifferential operators} in the sense of Convex Analysis. 

On page~199 in his monograph \cite{Si2} from 2008,
Stephen Simons asked the question whether or not $A+N_C$ is maximal
monotone when $A\colon \dom A\to X^*$ is linear and maximal monotone
and Rockafellar's constraint qualification $\dom A \cap \inte
C\neq\varnothing$ holds. 
In this manuscript, we provide an affirmative answer to Simons' question.
In fact, maximality of $A+N_C$ is guaranteed even when $A$ is a maximal
monotone linear relation, i.e., $A$ is simultaneously 
a maximal monotone operator and
a linear relation.

The paper is organized as follows.
In Section~\ref{s:aux}, we collect auxiliary results for future reference
and for the
reader's convenience. The main result (Theorem~\ref{t:main}) is proved
in Section~\ref{s:main}.

\section{Auxiliary Results} 
\label{s:aux}

\begin{fact}[Rockafellar] \label{f:F4}
\emph{(See {\cite[Theorem~3(a)]{Rock66}}, 
{\cite[Corollary~10.3]{Si2}}, or
{\cite[Theorem~2.8.7(iii)]{Zalinescu}}.)}\\
Let $f$ and $g$ be proper convex functions from $X$ to $\RX$.
Assume that there exists a point $x_0\in\dom f \cap \dom g$
such that $g$ is continuous at $x_0$. 
Then for every $z^*\in X^*$,
there exists $y^*\in X^*$ such that
\begin{equation}
(f+g)^*(z^*) = f^*(y^*)+g^*(z^*-y^*).
\end{equation}
\end{fact}

\begin{fact}[Fitzpatrick] 
\emph{(See {\cite[Corollary~3.9]{Fitz88}}.)}
\label{f:Fitz}
Let $A\colon X\To X^*$ be maximal monotone,  and set
\begin{equation}
F_A\colon X\times X^*\to\RX\colon
(x,x^*)\mapsto \sup_{(a,a^*)\in\gra A}
\big(\scal{x}{a^*}+\scal{a}{x^*}-\scal{a}{a^*}\big),
\end{equation}
which is the \emph{Fitzpatrick function} associated with $A$.
Then for every $(x,x^*)\in X\times X^*$, the inequality 
$\scal{x}{x^*}\leq F_A(x,x^*)$ is true,
and equality holds if and only if $(x,x^*)\in\gra A$. 
\end{fact}

\begin{fact}[Simons]\label{SL} 
\emph{(See {\cite[Corollary~28.2]{Si2}}.)}
Let $A\colon X\rightrightarrows X^*$ be maximal monotone. Then
\begin{equation}
\overline{\spand(P_X\dom F_A)}=\overline{\spand\dom A},
\end{equation}
where $P_X:X\times X^*\rightarrow X:(x,x^*)\mapsto x$.
\end{fact}

\begin{fact} \label{f:feb1}
Let $A\colon X\To X^*$ be a monotone linear relation, and set
\begin{equation}
(\forall x\in X)\quad 
q_A(x) = 
\begin{cases}
\thalb\scal{x}{Ax}, &\text{if $x\in\dom A$;}\\
\pinf, &\text{otherwise.}
\end{cases}
\end{equation}
Then $q_A$ is single-valued, convex, and nonnegative; in fact,
for $x$ and $y$ in $\dom A$, and $\lambda\in\RR$, we have 
\begin{align}
\lambda q_A(x) + (1-\lambda)q_A(y) - q_A(\lambda x + (1-\lambda)y) &=
\lambda(1-\lambda)q_A(x-y)\\
&= \thalb\lambda(1-\lambda)\scal{x-y}{Ax-Ay}.\notag
\end{align}
\end{fact}
\begin{proof}
This is a consequence of 
\cite[Proposition~2.2(iv) and Proposition~2.3]{BWY3}. 
While the results there 
are formulated in a reflexive Banach space, the proofs 
carry over \emph{verbatim} to the present general Banach space setting. 
\end{proof}

\begin{lemma}\label{l:l1} Let $C$ be  a nonempty closed convex
subset of $X$ such that $\inte C\neq \varnothing$.
Let $c_0\in \inte C$ and suppose that $z\in X\smallsetminus C$. 
Then there exists
$\lambda\in\left]0,1\right[$ such
that $\lambda c_0+(1-\lambda)z\in\bd C$.
\end{lemma}
\begin{proof}
Let $\lambda=\inf\big\{t\in[0,1]\mid tc_0+(1-t)z\in C\big\}$. Since $C$ is closed,
\begin{align} \lambda=\min\big\{t\in[0,1]\mid tc_0+(1-t)z\in C\big\}.\label{L:2}\end{align}
Because $z\notin C$, 
$\lambda> 0$. We now show that $\lambda c_0+(1-\lambda)z\in\bd C$.
Assume to the contrary that 
$\lambda c_0+(1-\lambda)z\in\inte C$. Then there exists 
$\delta\in\left]0,\lambda\right[$
such that $\lambda c_0+(1-\lambda)z-\delta(c_0-z)\in C$.
Hence
$(\lambda-\delta) c_0+(1-\lambda+\delta)z\in C$,
which contradicts \eqref{L:2}.
Therefore, 
$\lambda c_0+(1-\lambda)z\in\bd C$. 
Since $c_0\notin\bd C$, we also have $\lambda < 1$.
\end{proof}

The following useful result is a variant of \cite[Theorem~2.14]{BW}. 

\begin{lemma} \label{l:weird}
Let $A\colon X\To X^*$ be a set-valued operator, 
let $C$ be a nonempty closed convex subset of $X$,
and let $(z,z^*)\in X\times X^*$.
Set
\begin{equation}
I_C\colon X\To X^*\colon x\mapsto \begin{cases}
\{0\}, &\text{if $x\in C$;}\\
\varnothing, &\text{otherwise.}
\end{cases}
\end{equation}
Then $(z,z^*)$ is monotonically related to $\gra(A+N_C)$
if and only if
\begin{equation}
\text{
$(z,z^*)$ is monotonically related to $\gra(A+I_C)$
\quad and \quad $z \in \bigcap_{a \in\dom A\cap C}\big(a+T_C(a)\big)$,
}
\end{equation}
where $(\forall a\in C)$ 
$T_C(a) = \menge{x\in X}{\sup\scal{x}{N_C(a)}\leq 0}$.
\end{lemma}
\begin{proof}
``$\Rightarrow$'':
Since $\gra I_C\subseteq \gra N_C$, it follows 
that $\gra(A+I_C)\subseteq\gra(A+N_C)$; 
consequently, $(z,z^*)$ is monotonically related to $\gra(A+I_C)$. 
Now assume that $a\in\dom A \cap C$, 
and let $a^*\in Aa$. Then $(a,a^*+N_C(a))\subseteq \gra(A+N_C)$ and hence
$\scal{a-z}{a^*+N_C(a)-z^*}\geq 0$. This implies
$\pinf > \scal{a-z}{a^*-z^*}\geq \scal{z-a}{N_C(a)}$.
Since $N_C(a)$ is a cone, it follows that $\scal{z-a}{N_C(a)}\leq 0$ and
hence $z\in a+T_C(a)$. 
``$\Leftarrow$'':
Assume that $a\in\dom A\cap C$. 
Then $Aa = (A+I_C)a$, which yields $\scal{z-a}{Aa-z^*}\leq 0$, and also
$z-a\in T_C(a)$, i.e., $\scal{z-a}{N_C(a)}\leq 0$. 
Adding the last two inequalities, we obtain
$\scal{z-a}{Aa+N_C(a)-z^*}\leq 0$, i.e., $\scal{a-z}{(A+N_C)(a)-z^*}\geq
0$. 
\end{proof}

\section{Main Result}
\label{s:main}

\begin{theorem}\label{t:main}
Let $A:X\To X^*$ be a maximal monotone linear relation, 
let $C$ be a nonempty closed convex subset of $X$,
and suppose that $\dom A \cap \inte C\neq \varnothing$.
Then $A+N_C$ is maximal monotone.
\end{theorem}

\begin{proof}
Let $(z,z^*)\in X\times X^*$ and suppose that
\begin{equation} \label{e:zz*}
\text{
$(z,z^*)$ is monotonically related to $\gra(A+N_C)$. 
}
\end{equation}
It suffices to show that
\begin{equation} \label{e:ourgoal}
(z,z^*)\in\gra(A+N_C).
\end{equation}
We start by setting 
\begin{align} \label{e:defoff}
f\colon X\times X^*&\to \RX\\
(x,x^*)&\mapsto \scal{x-z}{x^*-z^*} + \iota_{\gra A}(x,x^*) + \iota_{C\times
X^*}(x,x^*)\notag\\
&= \big(\scal{x}{x^*}+\iota_{\gra A} + \iota_{C\times X^*}\big)
+\scal{(x,x^*)}{(-z^*,-z)} + \scal{z}{z^*}.\notag
\end{align}
If $(x,x^*)\in\dom f$, then $(x,x^*)\in\gra A$ and $x\in C$;
hence $x^*\in(A+N_C)x$ and thus $(x,x^*)\in\gra(A+N_C)$. 
In view of \eqref{e:zz*} and \eqref{e:defoff}, we deduce that 
$0 \leq \inf f(X\times X^*) = -f^*(0,0)$. Hence
\begin{equation} \label{e:f*neg}
f^*(0,0) \leq 0.
\end{equation}
Now let $q_A$ be as in Fact~\ref{f:feb1}. 
Since $\gra A$ is linear and hence convex,
it follows from Fact~\ref{f:feb1} that the function 
\begin{equation}\label{e:defofg}
g\colon X\times X^* \to \RX\colon
(x,x^*)\mapsto 2q_A(x) + \iota_{\gra A}(x,x^*) = \scal{x}{x^*} +
\iota_{\gra A}(x,x^*)
\end{equation}
is convex. 
Then
\begin{equation} \label{e:defofh}
h = g + \iota_{C\times X^*}
\end{equation}
is convex as well. 
Let
\begin{equation} \label{e:defofc0}
c_0 \in \dom A \cap \inte C,
\end{equation}
and let $c_0^*\in Ac_0$. 
Then $(c_0,c_0^*)\in \gra A \cap (\inte C \times X^*) = \dom g \cap
\intdom\iota_{C\times X^*}$, and $\iota_{C\times X^*}$ is continuous at
$(c_0,c_0^*)$. By Fact~\ref{f:F4} 
(applied to $g$ and $\iota_{C\times X^*}$),
there exists $(y^*,y^{**})\in X^{*}\times X^{**}$ such that
\begin{align} \label{e:puncha}
h^*(z^*,z) &= g^*(y^*,y^{**}) + \iota_{C\times X^*}^*(z^*-y^*,z-y^{**})\\
&=g^*(y^*,y^{**}) + \iota_{C}^*(z^*-y^*) + \iota_{\{0\}}(z-y^{**}).\notag
\end{align}
On the other hand, \eqref{e:defoff}, \eqref{e:defofg}, and \eqref{e:defofh}
imply that $h = f+\scal{\,\cdot\,}{(z^*,z)} -\scal{z}{z^*}$. Hence
$h^* = \scal{z}{z^*} + f^*(\,\cdot\, - (z^*,z))$, which, using \eqref{e:f*neg}, 
yields in particular 
\begin{equation} \label{e:punchb}
h^*(z^*,z) = \scal{z}{z^*} + f^*(0,0) \leq \scal{z}{z^*}.
\end{equation}
Combining \eqref{e:puncha} with \eqref{e:punchb}, we obtain
\begin{equation} \label{e:punchc}
g^*(y^*,y^{**}) + \iota_C^*(z^*-y^*) + \iota_{\{0\}}(z-y^{**}) \leq
\scal{z}{z^*}.
\end{equation}
Therefore, $y^{**}=z$ and
$g^*(y^*,z) + \iota_C^*(z^*-y^*) \leq \scal{z}{z^*}$.
Since $g^*(y^*,z) = F_A(z,y^*)$, we deduce that 
$F_A(z,y^*) + \iota_{C}^*(z^*-y^*)\leq \scal{z}{z^*}$;
equivalently,
\begin{equation} \label{e:punchd}
(\forall c\in C)\quad
F_A(z,y^*)-\scal{z}{y^*} + \scal{c-z}{z^*-y^*} \leq 0.
\end{equation}
We now claim that
\begin{equation}\label{e:smallgoal}
z\in C.
\end{equation}
Assume to the contrary that \eqref{e:smallgoal} fails, i.e., that 
$z\notin C$.
By \eqref{e:punchd}, $(z,y^*)\in\dom F_A$. 
Using Fact~\ref{SL} and the fact that $\dom A$ is a linear subspace of $X$,
we see that
$z\in P_X(\dom F_A) \subseteq \overline{\spand P_X(\dom F_A)}
= \overline{\spand\dom A} = \overline{\dom A}$. 
Hence there exists a sequence $(z_n)_\nnn$ in 
$(\dom A)\smallsetminus C$ such that 
$z_n\to z$. 
By Lemma~\ref{l:l1}, 
$(\forall\nnn)$ $(\exists \lambda_n\in\left]0,1\right[\,)$ 
$\lambda_n z_n+(1-\lambda_n)c_0\in\bd C$.
Thus,
\begin{equation}
\label{LZZ:0}
(\forall\nnn)\quad 
\lambda_n z_n+(1-\lambda_n)c_0\in\dom A \cap \bd C.
\end{equation}
After passing to a subsequence and relabeling if necessary,
we assume that $\lambda_n\to\lambda\in[0,1]$.
Taking the limit in \eqref{LZZ:0}, we deduce that
$\lambda z + (1-\lambda)c_0\in\bd C$.
Since $c_0\in\inte C$ and $z\in X\smallsetminus C$,
we have $0<\lambda$ and $\lambda<1$. 
Hence 
\begin{equation} \label{e:zeroun}
\lambda_n \to \lambda \in\zeroun.
\end{equation}
Since $\inte C\neq\varnothing$, 
Mazur's Separation Theorem (see, e.g., \cite[Theorem~2.2.19]{Megg}) 
yields a sequence $(c^*_n)_\nnn$ in $X^*$ 
such that
\begin{equation} \label{e:punchf}
(\forall\nnn)\quad
c_n^*\in N_C\big(\lambda_n z_n+(1-\lambda_n)c_0\big)
\;\;\text{and}\;\;
\|c_n^*\|_* = 1.
\end{equation}
Since $c_0\in\inte C$, there exists $\delta>0$ such that
$c_0+\delta B_X\subseteq C$. 
It follows that
\begin{equation} \label{e:punchg}
(\forall\nnn) \quad
\delta \leq \lambda_n\scal{z_n-c_0}{c_n^*}.
\end{equation}
Since the sequence $(c_n^*)_\nnn$ is bounded,
we pass to a weak* convergent \emph{subnet} 
$(c^*_\gamma)_{\gamma\in\Gamma}$, say
$c^*_\gamma \weakstarly c^*\in X^*$.
Passing to the limit in \eqref{e:punchg} along subnets, we see that
$\delta \leq \lambda\scal{z-c_0}{c^*}$; hence, using \eqref{e:zeroun}, 
\begin{equation} \label{e:punchh}
0 < \scal{z-c_0}{c^*}.
\end{equation}
On the other hand and borrowing the notation of Lemma~\ref{l:weird},
we deduce from \eqref{LZZ:0}, \eqref{e:zz*}, 
and Lemma~\ref{l:weird} that 
$(\forall\nnn)$ $z\in (\Id+T_C)(\lambda_nz_n+(1-\lambda_n)c_0)$, 
which in view of \eqref{e:punchf} yields
\begin{equation} \label{e:punchi}
(\forall\nnn)\quad
\scal{z-(\lambda_nz_n+(1-\lambda_n)c_0)}{c_n^*} \leq 0.
\end{equation}
Taking limits in \eqref{e:punchi} along subnets, we deduce
$\scal{z-(\lambda z+(1-\lambda)c_0)}{c^*}\leq 0$. 
Dividing by $1-\lambda$ and recalling \eqref{e:zeroun}, we thus have
\begin{equation} \label{e:punchj2}
\scal{z-c_0}{c^*}\leq 0.
\end{equation}
Considered together, the inequalities 
\eqref{e:punchh} and \eqref{e:punchj2} are absurd --- 
we have thus verified \eqref{e:smallgoal}. 

Substituting \eqref{e:smallgoal} into \eqref{e:punchd}, 
we deduce that
\begin{equation}
F_A(z,y^*)\leq\scal{z}{y^*}.
\end{equation}
By Fact~\ref{f:Fitz}, 
\begin{equation}
\label{e:punchx}
(z,y^*)\in\gra A
\end{equation}
and $F_A(z,y^*)=\scal{z}{y^*}$.
Thus, using \eqref{e:punchd} again, we see that
$\sup_{c\in C} \scal{c-z}{z^*-y^*}\leq 0$, i.e., that
\begin{equation}
\label{e:punchy}
(z,z^*-y^*)\in \gra N_C. 
\end{equation}
Adding \eqref{e:punchx} and \eqref{e:punchy}, we obtain
\eqref{e:ourgoal}, and this completes the proof.
\end{proof}

\begin{corollary}\label{c:main}
Let $A:X\To X^*$ be maximal monotone and at most single-valued,
and let $C$ be a nonempty closed convex subset of $X$.
Suppose that $A|_{\dom A}$ is linear, 
and that $\dom A \cap \inte C\neq \varnothing$.
Then $A+N_C$ is maximal monotone.
\end{corollary}

\begin{remark}
Corollary~\ref{c:main} provides an affirmative answer
to a question Stephen Simons raised in his 2008 monograph
\cite[page~199]{Si2} concerning \cite[Theorem~41.6]{Si}. 
\end{remark}

\section*{Acknowledgment}
Heinz Bauschke was partially supported by the Natural Sciences and
Engineering Research Council of Canada and
by the Canada Research Chair Program.
Xianfu Wang was partially supported by the Natural
Sciences and Engineering Research Council of Canada.


\end{document}